\documentclass[a4paper, 12pt]{amsart}

\usepackage {algorithmic, algorithm, epsfig,nicefrac,amsthm, amssymb}

\usepackage[utf8]{inputenc}

\usepackage{mathptmx}

\newtheorem{lemma}{Lemma}[section] \newtheorem{corollary}[lemma]{Corollary} \newtheorem{theorem}[lemma]{Theorem}

\theoremstyle{definition} \newtheorem{definition}[lemma]{Definition}

\def\ZZ{{\mathbb Z}}
\def\RR{{\mathbb R}}

\DeclareMathOperator{\rem}{rem}

\DeclareMathOperator{\Hilb}{Hilb}
\DeclareMathOperator{\ld}{ld}

\DeclareMathOperator{\parf}{par}

\DeclareMathOperator{\fFunc}{f}
\DeclareMathOperator{\gFunc}{g}

\newcommand{\bigO}{\mathcal{O}}

\def\ln{\operatorname{log}}

\DeclareRobustCommand{\bigO}{%
  \text{\usefont{OMS}{cmsy}{m}{n}O}%
}

\begin {document}

\title[Polynomial-size vectors for unimodular triangulation]{Polynomial-size vectors are enough for the unimodular triangulation of simplicial cones}

\author{Michael von Thaden}
\address{FH Westk\"{u}ste,
	25746 Heide,
	Germany}
\email{\texttt{v\_thaden@t-online.de}}



\maketitle

\begin{abstract}

In a recent paper, Bruns and von Thaden established a bound for the length of vectors involved in a unimodular triangulation of simplicial cones. The bound is exponential in the square of the logarithm of the multiplicity, and improves previous bounds significantly. In this paper we will prove that a bound, which is polynomial in the multiplicity $\mu$, exists. In detail, the bound is of the type $\mu^{\fFunc(d)}$ with $\fFunc(d) \in \bigO(d)$. 

\end{abstract}

%\fontsize{16}{19}\normalfont
\marginparwidth=5em
\marginparpush=10pt

\section{Introduction}

Unimodular triangulations of polytopes, cones and polyhedral complexes are a useful and important tool in many subfields of mathematics like algebraic geometry, commutative algebra, (enumerative) combinatorics or integer programming. In toric geometry, unimodular triangulations of cones correspond to desingularizations of a toric varieties. Here, the standard method for desingularization normally leads to triangulations which involve rather long vectors.  

In \cite{BvT17} Bruns and von Thaden established a bound for the length of vectors involved in a unimodular triangulation of simplicial cones. Length is hereby measured by the basic simplex $\Delta_C$ of $C$ that is spanned by the origin and the extreme integral generators of $C$. We are interested in an upper bound for the dilatation factor $c$ for which all subdividing vectors are contained in $c\Delta_C$. Bruns and von Thaden gave an upper bound for $c$, which was exponential in the square of the logarithm of the multiplicity of the cone $C$, hereby improving a result from Bruns and Gubeladze \cite[Theorem 4.1]{BG02a} which itself was a slight improvement of the standard argument applied for the desingularization of toric varieties. Bruns and von Thaden mentioned in \cite{BvT17} that the next goal would be a bound that is polynomial in the multiplicity. In this paper we will prove that such a bound, which is polynomial in the multiplicity $\mu$, indeed exists: the bound is of type $\mu^{\fFunc(d)}$ with $\fFunc(d) \in \bigO(d)$.

Of course, a corresponding result for the unimodular triangulation of lattice polytopes would be very desirable but this seems currently out of reach. The best result so far is the celebrated Knudsen-Mumford-Waterman theorem \cite{KKMS73}. It states that a $c'$ exists such that the multiples $cP$ of a lattice polytope $P$ have unimodular triangulations for \emph{all} $c\geq c'$. But Knudsen, Mumford and Waterman did not provide an explicit bound. Recently, Haase, Paffenholz, Piechnik and Santos closed this gap in \cite{HPPS} and provided an explicit bound, which is doubly exponential in the volume of the lattice polytope $P$. Interestingly, if one is only interested in unimodular covers of lattice polytopes instead of unimodular triangulations one can do much better. In \cite{BG09} Bruns and Gubeladze showed that multiples $cP$ of lattice polytopes $P$ can be covered by unimodular simplices for all $c\geq \gFunc(d)$ with $\gFunc \in \bigO(d^6)$. So, this threshold does only depend on $d$ and not on the multiplicity of $P$.

\cite{HPPS} gives a comprehensive overview of the topic of unimodular triangulations as does \cite[Chapter 3]{BG09}. Furthermore, we refer the reader to \cite{BG09} for any unexplained terminology.

\section{Auxiliary results}

One of the main ideas of the proof in \cite{BvT17} was that a cone whose multiplicity is a power of 2, or, generally speaking, whose multiplicity is a product of small primes could be triangulated using just short vectors. 
Therefore, in the first step one might wish to triangulate a cone into subcones whose multiplicities are exclusively products of small primes while keeping the subdividing vectors as short as possible.

If one wants to apply stellar subdivision to come up with a triangulation of $C$ by cones of the desired type, what kind of vectors should be used in the stellar subdivisions?
Recall that if the primitive vectors $v_1,\dots,v_d\in\ZZ^d$ generate a simplicial cone $C$ of dimension $d$, and if $U$ denotes the sublattice of $\ZZ^d$ spanned by these vectors, then $\mu(C)$ is the index of $U$ in $\ZZ^d$, and each residue class has a representative in 
$$
\parf(v_1,\dots,v_d)=\{q_1v_1+\dots+q_dv_d: 0\le q_i < 1\}.
$$
If $p$ divides $\mu(C)$, then there is an element of order $p$ in $\ZZ^d/U$, and consequently there exists a vector
\begin{equation}\label{vectorform}
x= \frac{1}{p}\sum_{i=1}^d z_iv_i \in \ZZ^d\setminus{\{0\}}, \qquad z_i\in \ZZ, \ 0\le z_i < p.
\end{equation}

If we do now apply stellar subdivision with respect to $x$ to $C$, then the resulting cones $E_i$ have multiplicities $\mu(E_i) = z_i/p\cdot \mu(C) <\mu(C)$. So, in essence, one substitutes a prime factor $p$ in the factorization of $\mu(C)$ by the number $z_i$. This means that if one could choose $x$ in a way that the $z_i$ are composite numbers, one increases the number of prime factors for the triangulating cones hereby ensuring that the prime factors are getting smaller.

But are there always vectors of type \eqref{vectorform} such that all the $z_i$ are composite numbers? In general, this is not always the case: let $d+1$ be a prime number and let the cone $C\subset \mathbb{R}^d$ be generated by the vectors 
$$
v_1= \sum_{j=2}^{d} (d+1-j) e_j +(d+1)e_1, \quad v_i = e_i, \quad i=2, \ldots, d.
$$ 
Let $\mathbb{P}$ be the set of all primes. Then $\mu(C) = d+1\in \mathbb{P}$. Furthermore, each residue class of $\ZZ^d$ modulo the sublattice $U$ generated by the vectors $v_i$ has a representative, which is of the form 
$$
x_j = \frac{1}{d+1} \sum_{k=1}^d (jk \ \rem d+1) v_k, \quad j= 1, \ldots, d+1, 
$$
where $a \ \rem b$ denotes the remainder of $a$ modulo $b$. Then, for every $j$ we have \linebreak[4] $\{jk \ \rem d+1 : k=1, \ldots, d\} = \{1, \ldots, d\}$. So, in this case there is no vector $x$ of the form \eqref{vectorform} such that all $z_i$ are just composite numbers.

But is there a condition which ensures that such a vector $x$ of type \eqref{vectorform} with all $z_i$ being composite numbers exists? 
We will now prove quite easily that 
$\parf(v_1,\dots, v_d)$ indeed always contains a vector of form \eqref{vectorform} such that all the $z_i$ are composite numbers -- as long as the largest prime factor  of $\mu(C)$ is bounded below by $e^{\tau d}$, where $\tau = 1.25506$. This fact is a direct consequence of the following lemma which has already been proved in \cite{BvT17} with the help of an upper bound for the prime number counting function $\pi(n)$, as provided by Rosser \& Schoenfeld in \cite{RS62}.

\begin{lemma}\label{avoid}
With the notation introduced, let $M\subset \{1,\dots,d\}$ such that 
$$
|M| \leq \frac{\ln(p)}{\tau},\quad \tau=1.25506.
$$
Then there exists an element $x$ of order $p$ modulo $U$ such that none of the coefficients $z_i$, $i\in M$, is an odd prime $<p$. 
\end{lemma}

If one takes $M = \{1, \ldots, d\}$, the lemma  implies that there exists an element $x \in \parf(v_1, \\ \ldots, v_d) \setminus\{0\}$ of type \eqref{vectorform} such that none of the coefficients $z_i, i=1, \ldots,d$ in \eqref{vectorform} is an odd prime as long as 
$$
p \geq e^{\tau d}, \quad \tau=1.25506.
$$
Therefore, we have
\begin{theorem}\label{avoidspecial}
Let $C = \RR_+ v_1 + \cdots + \RR_+ v_d \subset \RR^d , \, d\geq 2$ be a simplicial $d$-cone such that 
$$
p_{\max}:=p_{\max}(\mu(C)) \geq e^{\tau d},
$$
where
$p_{\max}(n):= \max\{p\in \mathbb{P}:p\mid n\}$.
Then there exists a vector 
$$
x = \frac{1}{p_{\max}}\left(\sum_{i=1}^d z_i v_i\right) \in \parf(v_1, \ldots, v_d) \setminus\{0\}.
$$
such that $z_i \notin \mathbb{P}_{>2}$ for all $i$.
\end{theorem}

Hence, as long as $\mu(C)$ has a prime factor $p \geq e^{\tau d}$, there is also an element $x = \frac{1}{p}\left(\sum_{i=1}^d z_i v_i\right) \in \parf(v_1, \ldots, v_d) \setminus\{0\}$  such that all $z_i$ are composite numbers or are equal to 2. 
These short vectors can then be used for successive stellar subdivision until one arrives at cones $D_i$ for which $p_{\max} (\mu(D_i)) < e^{\tau d}$ and which do constitute a triangulation of the original cone $C$. 

The following definition will help us to shorten any further explanations or statements for this kind of triangulation procedures.

\begin{definition}
An $f$-triangulation is defined as a triangulation of a cone $C$ by cones $D_i$ for which $p_{\max}(\mu(D_i)) <f$ for all $i$.
\end{definition}

\section{The algorithm}

With the help of Theorem \ref{avoidspecial} we are now ready to formulate an algorithm which provides us with an $e^{\tau d}$-triangulation of a cone $C$ by cones $D_i$. As we will see, the vectors involved in this triangulation are short and the multiplicities of the cones $D_i$ in the resulting triangulation are smaller than the multiplicity $\mu(C)$ of the original cone $C$.
\pagebreak[2]

\noindent \hrulefill

\noindent\textbf{Bounded prime factors triangulation -- BPFT}
\vspace{-3 mm}

\noindent \hrulefill

\begin{algorithmic}[1]

\REQUIRE The initial cone $C = \RR_+ v_1 + \cdots + \RR_+ v_d \subset \RR^d$

\ENSURE An $e^{\tau d}$-triangulation $\hat{T}(C)$ of $C$

\STATE $\hat{T}(C):= \{C\}$

\STATE $\hat{A}(C):=\{C\}$

\STATE $\xi_C(-i) := v_i$ for $i =1, \ldots , d$

\STATE $\xi_C(i) := 0$ for $i \in \mathbb{N}_0$

\WHILE {$\hat{T}(C)$ contains a cone $D= \RR_+ \xi_D(i_1)+ \cdots + \RR_+ \xi_D(i_d)$ (where $i_1 >i_2> \ldots >i_d\geq -d$) such that $p_{\max}(\mu(D)) \geq e^{\tau d}$}

\STATE $p := p_{\max}(\mu(D))$

\STATE {\textbf{FIND} $x = \nicefrac{1}{p}\left(\sum_{j=1}^d z_j \xi_D(i_j)\right)\in \parf(\xi_D(i_1), \ldots, \xi_D(i_d))\setminus\{0\}$ such that $z_j \notin \mathbb{P}_{>2}$ for all $j$ (exists due to Theorem \ref{avoidspecial})}\\

\FORALL {$E \in \hat{T}(C)$ with $x \in E$}

\STATE Apply stellar subdivision to $E$ by $x$ (let
$E_j$ $(j= 1, \ldots, m)$ be the resulting cones)

\STATE $\hat{T}(C):= (\hat{T}(C) \setminus \{E\}) \cup \{E_j : \, j= 1, \ldots, m\}$

\STATE $\hat{A}(C):= \hat{A}(C) \cup \{E_j : \, j= 1, \ldots, m\}$

\ENDFOR

\STATE $\nu:= \max\{i : \xi_E (i)\neq 0\}$

\FORALL{$j = 1, \ldots, m$}

\FORALL{$k \leq \nu$}

\STATE $\xi_{E_j}(k) := \xi_{E}(k)$

\ENDFOR

\STATE $\xi_{E_j}(\nu +1) :=x$

\ENDFOR

\ENDWHILE

\STATE Return $\hat{T}(C)$
\end{algorithmic}

\vspace{-3 mm}
\noindent \hrulefill
\vspace{5 mm}

For a simplicial $d$-cone $C$ the BPFT algorithm computes an $e^{\tau d}$-triangulation of $C$. It applies successive stellar subdivisions to the initial cone $C$ and it stops when all multiplicities only have prime factors smaller than $e^{\tau d}$. Finally, it stops after finitely many iterations, because, if $E$ results from $D$ by stellar subdivision in the course of the BPFT algorithm, then we have $\mu(E) < \mu(D)$.

As in \cite{BvT17} the set $\hat{A}(C)$ contains the original cone $C$ and all cones being created in the course of the algorithm and the set $\hat{T}(C)$ is a strict subset of $\hat{A}(C)$ unless $\mu(C)$ is not divisible by a prime greater than or equal to $e^{\tau d}$. $\hat{A}(C)$ has been introduced out of technical reasons; it will help us to analyze certain properties of the resulting triangulation. The vectors $\xi_D(i)$ for a cone $D$ include all extremal generators of all cones $E$ containing the cone $D$. In particular, they also include the extremal generators of the cone $D$ itself. 
 
In section \ref{BoundsETAU} we will show that the generators of the resulting cones $E_i\in\hat{T}(C)$ are short. Building on these results we will finally, in section \ref{BoundsUT}, introduce new bounds for the length of vectors involved in unimodular triangulations of simplicial cones.

\section{Bounds for $e^{\tau d}$-triangulation} \label{BoundsETAU}

\begin{theorem}\label{vectorlength1}
Let $D \in \hat{T}(C)$. Then, for all $s\geq 0$
$$
\xi_D(s) \in \biggl(d \cdot 2^s \biggr)\Delta_C.
$$
\end{theorem}

\begin{proof}
The proof of this theorem is similar to the proof of Theorem 4.1 in \cite{BG02a} and the proof of 
Theorem 2.1 in \cite{BvT17}.
We consider the following sequence:
$$
h_k = 1, \quad k \leq -1,\qquad
h_k = h_{k-1} + \cdots + h_{k-d}, \quad k\geq 0.
$$
\noindent Because
$$
h_k - h_{k-1} = h_{k -1} -  h_{k-d-1}
$$
for $k\geq 1$ and
$h_0 \geq h_l $ for $l \leq -1$, it follows by induction that this sequence is increasing.
Since for $k\geq 1$
$$
h_k = h_{k-1} + ( h_{k-2} + \cdots + h_{k-d-1}) -  h_{k-d-1} =
2h_{k-1} - h_{k-d-1} < 2h_{k-1},
$$
and because $h_0 =d$, 
we arrive at
$$
h_k \leq d \cdot 2^k
$$
for $k \geq 0$. This inequality will be needed in the following.

Now, we will prove via induction on $s$ that
$$
\xi_D(s) \in h_s\Delta_C
$$
So, let $s =0$. If $\xi_D(0) = 0$, there is nothing to prove. 

So, suppose that $\xi_D(0) \neq 0$. By the construction of $\xi_D(0)$ it follows that this vector was used for the stellar subdivision of the initial cone $C$. Hence, $\xi_D(0)$ is of the form
$$
\xi_D(0) = \sum_{i = 1}^d z_i v_i \in \ZZ^d \setminus \{0\}.
$$
where
$z_i < 1$ for all $i$. 
Therefore, $x \in d\Delta_C$, which finishes the case $s=0$. 

For the induction step assume the statement is true for $s-1\geq 0$. Again there is nothing to prove if $\xi_D(s) = 0$. Otherwise $\xi_D(s) \neq 0$ is a vector used for stellar subdivision. With the same notation as above, it follows
by construction of $\xi_D(s)$, that
$$
\xi_D(s) = \sum_{i = 1}^d z_i \xi_D(j_i) \in \ZZ^d \setminus \{0\}
$$
such that $s> j_1>j_2> \ldots > j_d$ and again $z_i < 1$. So, it follows by induction that
$$
\xi_D(s) \in (h_{j_1}+ \cdots + h_{j_d}) \Delta_C.
$$
Because the $h_i$ are increasing, this means that
$$
\xi_D(s) \in (h_{s-1}+ \cdots + h_{s-d}) \Delta_C = h_s \Delta_C,
$$
which finishes the proof.
\end{proof}

The next definition will be helpful in showing that the length of every chain of cones
$$
E_0=D \subset E_1 \subset E_2 \ldots \subset E_L=C,
$$ 
where $E_i$ is generated from $E_{i+1}$ by stellar subdivision and $D$ belongs to the resulting $e^{\tau d}$-triangulation of $C$, is relatively short.

\begin{definition}\label{Defphi1}
Let $n$ be a natural number, $n= \prod_{i=1}^{\infty} p_i^{\alpha_i}$ be its prime
decomposition.
Then  we define $\phi(n)= \ld(n)- \eta(n)$, where $\eta(n) = \sum_{i=1}^{\infty} \alpha_i$. Hence, $\phi(n) = \sum_{i=1}^{\infty}\alpha_i\left(\ld(p_i)- 1 \right)$.
\end{definition}

The function $\phi$ has some obvious nice properties.

\begin{lemma}\label{phi1}\leavevmode
\begin{enumerate}
\item $\phi(ab) = \phi(a) + \phi(b)$ for $a,b \in \mathbb{N}$,
\item $\phi(a/b) = \phi(a) - \phi(b)$ for $a,b \in \mathbb{N}$, $b\mid a$,
\item $\phi(a) \geq 0$ for $a \in \mathbb{N}$.
\end{enumerate}
\end{lemma}

\begin{lemma}\label{phimu}
Let $D,E \in \hat{A}(C)$ such that $E$ results from $D$ by stellar subdivision in the course of the BPFT algorithm. Then
$$
\phi(\mu(E)) \leq \phi(\mu(D))-1.
$$
\end{lemma}

\begin{proof}
Due to lines 7 and 9 of the algorithm,
$$
\mu(E) = \mu(D)\frac{f}{p_{\max}},
$$
where $p_{\max}=\max\{p\in \mathbb{P}:p|\mu(D)\}$ and $f\in\mathbb{N}_{>0}$ is either 
\begin {enumerate}
\item a composite number smaller than $p_{\max}$, i.e. $f= uv$ with natural numbers $u, v >1$ or
\item $f=2 < p_{\max}$.
\end{enumerate}

For the first case we have by Lemma \ref{phi1} and because $p_{\max}\mid\mu(D)$
$$
\phi(\mu(E)) = \phi(\mu(D))-\phi(p_{\max}) + \phi(f)= 
$$
$$
\phi(\mu(D)) + \phi(u)+\phi(v)-\ld(p_{\max})-1\leq \phi(\mu(D)) - 1.
$$

For the second case it follows that
$$
\phi(\mu(E)) = \phi(\mu(D))-\phi(p_{\max}) + \phi(2)= 
$$
$$
\phi(\mu(D)) - \phi(p_{\max}) \leq \phi(\mu(D)) - 1,
$$
because $p_{\max} \geq e^{\tau d}\geq 3.5$ for all $d\geq 1$, which implies that  $p_{\max} \geq 5$, because $p_{\max}$ is a prime. Therefore, $\phi(p_{\max}) \geq \phi(5) = \ld(5)-1>1$. This proves the lemma.
\end{proof}

\begin{lemma}\label{chicone1}
Let $D\in \hat{A}(C)$ be an arbitrary cone resulting from the BPFT algorithm. Furthermore, we define
$$
\chi(D)= \max \{i: \xi_D(i) \neq 0\}.
$$
Then
$$
\chi(D) \leq \phi(\mu(C))-1.
$$
\end{lemma}

\begin{proof}
Let $D\in \hat{A}(C)$. By the algorithm, there is chain of cones
$$
E_0=D \subset E_1 \subset E_2 \ldots \subset E_L=C
$$ 
such that $E_i$ is generated from $E_{i+1}$ by stellar subdivision. Lemma \ref{phimu} implies that $\phi(\mu(D)) \leq \phi(\mu(C))-L$. On the other hand, by construction, $\chi(D) = \chi(C)+L$, where $\chi(C) =-1$. Therefore
$$
\chi(D) =L-1\leq \phi(\mu(C))-\phi(\mu(D))-1.
$$ 
This proves the lemma, because $\phi(a)\geq 0$ for all $a \in \mathbb{N}$.
\end{proof}
 
\begin{corollary} \label{CorrBPFT}
Every simplicial $d$-cone
$C = \RR_+ v_1 + \cdots + \RR_+ v_d \subset \RR^d , \, d\geq 2$, which is not already unimodular (i.e., $\mu(C)>1$) has an $e^{\tau d}$-triangulation $C = D_1 \cup \ldots \cup D_t$ such that 
$$
\Hilb(D_i) \subset \left(\frac{d}{4} \cdot \mu(C)\right)\Delta_C
$$ 
for all $i$.
\end{corollary}
\begin{proof}
Due to \ref{vectorlength1} and \ref{chicone1} it follows that
$$
\Hilb(D_i) \subset \left(d\cdot2^{\phi(\mu(C))-1}\right)\Delta_C.
$$
Because due to \ref{Defphi1} and $\mu(C) >1$ we have $\phi(\mu(C)) \leq \ld(\mu(C)) - 1$. It follows that $2^{\phi(\mu(C))-1} \leq \mu(C)/4$, which finally proves the corollary.
\end{proof}

\section{Bounds for unimodular triangulation}\label{BoundsUT}

Building on the previous  bound we will now introduce new bounds for the length of vectors involved in unimodular triangulations of simplicial cones. This will be done with the help of the following corollary from \cite{BvT17}.

\begin{theorem}\label{simpler}
Let $\epsilon= 5+3/2\cdot \ld(3/2)$ and $\rho = 1/2 \cdot \ld(3/2)$. So, $\epsilon \approx 5.88$ and $\rho \approx 0.29$.
Then every simplicial $d$-cone
$C = \RR_+ v_1 + \cdots + \RR_+ v_d \subset \RR^d , \, d\geq 2$, which is not already unimodular (i.e., $\mu(C)>1$) has a unimodular triangulation $C = D_1 \cup \ldots \cup D_t$ such that for all $i$
$$
\Hilb(D_i) \subset \left( \frac{d^2}{64} \cdot \mu(C)^{\rho \cdot \ld(\mu(C)) + \epsilon} \right)
\Delta_C.
$$
\end{theorem}

Furthermore, we will need the following lemma, which will help us with connecting the previous corollary and Corollary \ref{CorrBPFT} to achieve our main result of a new upper bound for the length of vectors involved in the unimodular triangulation of simplicial cones.

\begin{lemma}\label{PrimeTransfer}
Let us assume that every simplicial $d$-cone $E$ with $\mu(E) = p$ admits a unimodular triangulation $E = F_1 \cup \ldots \cup F_t$ such that 
$$
\Hilb(F_i) \subset k_{p,d}\Delta_E.
$$
for all $i$ for a certain $k_{p,d} \in \mathbb{R}$.
Let $C = \RR_+ v_1 + \cdots + \RR_+ v_d \subset \RR^d , \, d\geq 2$, be a simplicial $d$-cone and let $p\in \mathbb{P}$ such that $p|\mu(C)$. Then $C$ admits a triangulation $C = D_1 \cup \ldots \cup D_t$ with $D_i = \RR_+ w_1^i + \cdots + \RR_+ w_d^i$ such that 
\begin{enumerate}
\item $w_j^i \in k_{p,d}\Delta_C$ and
\item $\mu(D_i) = \mu(C)/p$
\end{enumerate}
for all $i$ and $j$.
\end{lemma}

\begin{proof}
Because $p|\mu(C)$, we know that there exists a vector
$$
x= \frac{1}{p}\sum_{i=1}^d z_iv_i \in \ZZ^d\setminus{\{0\}}, \qquad z_i\in \ZZ, \ 0\le z_i < p.
$$
W.l.o.g. we can assume that $z_1 = 1$. Now, let $E=\RR_+ v'_1 + \cdots + \RR_+ v'_d \subset \RR^d , \ d\geq 2$, be the simplicial $d$-cone generated by the vectors 
$$
v_1'= \sum_{j=2}^{d} (p-z_j) e_j + p e_1, \quad v_i' = e_i, \quad i=2, \ldots, d.
$$
Then $\mu(E) = p$. Furthermore, 
$$
x'= \frac{1}{p}\sum_{i=1}^d z_i v'_i = \sum_{i=1}^{d} e_i \in \ZZ^d\setminus{\{0\}}.
$$
Let  $F_i = \RR_+ {w^i_1}' + \cdots + \RR_+ {w^i_d}' \subset \RR^d$ be cones which constitute a unimodular triangulation of $E = F_1 \cup \ldots \cup F_t$ such that $\Hilb(F_i) \subset k_{p,d} \Delta_E$ for all $i$. Because $\mu(E)=p$ is the index of the sublattice $U$, which is spanned by the vectors $v_1', \ldots, v_d'$ , in $\mathbb{Z}^d$, it follows that $\ZZ^d$ modulo $U$ is generated by each non-null element. One representative of such an element is obviously $x'$. So, since ${w^i_j}' \in E  \cap \ZZ^d$,  it follows that for all $i$ and $j$ there exists a $l^i_j \in \mathbb{N}$ such that 
$$
{w^i_j}'  =   \frac{1}{p}\left(\sum_{k=1}^d a^i_{jk} v'_k\right),
$$ 
where $a^i_{jk}\in \mathbb{N}$ and
$$
(a^i_{jk}-l^i_jz_k) \ \rem p = 0. 
$$ 
for all $k$. This implies that also 
$$
w^i_j  =   \frac{1}{p}\left(\sum_{k=1}^d a^i_{jk} v_k\right) \in \ZZ^d
$$
for all $i$ and $j$.

Furthermore, because we have $\Hilb(F_i) \subset k_{p,d} \Delta_E$ and because the $F_i$ are unimodular, it follows that ${w^i_j}' \in k_{p,d} \Delta_E$ for all $i$ and $j$. Therefore, we also have  
$$
w^i_j \in k_{p,d}\Delta_{C}.
$$

Let now ${W^i}' \in \RR^{d\times d}$ be the matrix formed by the row vectors ${w^i_j}'$, let \linebreak[4] $A^i:=(a^i_{jk})_{j=1, \ldots, d, k= 1, \ldots, d}$ and $V' \in \RR^{d\times d}$  be the matrix formed by the row vectors $v'_k$. Then we have that 
$$
\det\left(\frac{1}{p} \cdot A^i V'\right) = \det\left({W^i}'\right) = \mu(F_i) = 1
$$ 
for all $i$, which implies that 
$$
\det \left(\frac{1}{p} \cdot A^i\right) = \frac{1}{\det(V')} = \frac{1}{\mu(E)} = \frac{1}{p}.
$$

Therefore, the triangulation of $C = D_1 \cup \ldots \cup D_t$ given by $D_i = \RR_+ w^i_1 + \cdots + \RR_+ w^i_d \subset \RR^d$ 
has the desired properties, because, first, we have that for $V \in \RR^{d\times d}$ formed by the row vectors $v_k$
$$
\mu(D_i) = \det\left(\frac{1}{p}\cdot A^i V\right) = \det\left(\frac{1}{p} \cdot A^i\right) \cdot \mu(C) = \frac{\mu(C)}{p}
$$ 
And second, we have already shown that $w_j^i \in k_{p,d}\Delta_C$ for all $i$ and $j$, if $\Hilb(F_i)\subset k_{p,d}\Delta_E$.
\end{proof}

\begin{corollary}\label{FinResI}
Let $\gamma = \rho\tau\ld(e)$ and $\kappa = \epsilon -5$. So, $\gamma \approx 0.53$ and $\kappa \approx 0.88$. Then every simplicial $d$-cone
$C = \RR_+ v_1 + \cdots + \RR_+ v_d \subset \RR^d , \, d\geq 2$, which is not already unimodular (i.e., $\mu(C)>1$) has a unimodular triangulation $C = D_1 \cup \ldots \cup D_t$ such that for all $i$
$$
\Hilb(D_i) \subset \left(\frac{d}{4} \cdot \mu(C)^{\gamma d+2\ld(d)+\kappa}\right)\Delta_C.
$$
\end{corollary}

\begin{proof}
Due to \ref{CorrBPFT} $C$ has an $e^{\tau d}$-triangulation $C = D_1 \cup \ldots \cup D_t$ such that for all $i$
$$
\Hilb(D_i) \subset \left(\frac{d}{4} \cdot\mu(C)\right)\Delta_C.
$$ 
Furthermore,
\begin{equation}\label{UpBoundMu}
\mu(D_i) \leq \mu(C)
\end{equation}
for all $i$.

So, let $\mu(D_i) = \prod_{j=1}^{n_i} p_j^{\alpha_{j,i}}$ be the prime decomposition of $\mu(D_i)$, where $p_1<\ldots <p_{n_i} < e^{\tau d}$. 
Then, due to successive application of Lemma \ref{PrimeTransfer}
and Corollary \ref{simpler} it follows that each of the cones $D_i$ admits a unimodular triangulation $D_i = F_1^i \cup \ldots \cup F_{s_i}^i$ such that
$$
\Hilb(F_k^i) \subset \prod_{j=1}^{n_i} \left( \frac{d^2}{64}p_j^{\rho\ld(p_j)+\epsilon} \right)^{\alpha_{j,i}}\Delta_{D_i}
$$
for all $i,k$.
Since the $D_i$ constitute an $e^{\tau d}$-triangulation of $C$, the latter has a unimodular triangulation 
$$
C = \bigcup_{i=1}^t \bigcup_{k=1}^{s_i} F_k^i
$$ 
such that, for all $i$ and $k$, we have
$$
\Hilb(F_k^i) \subset \left(\frac{d}{4} \cdot\mu(C)\cdot \prod_{j=1}^{n_i} \left( \frac{d^2}{64}p_j^{\rho\ld(p_j)+\epsilon} \right)^{\alpha_{j,i}}\right)\Delta_{C}.
$$

Because
$$
\sum_{j=1}^{n_i} \alpha_{j,i} \leq \ld(\mu(D_i))
$$
it follows that
$$
\prod_{j=1}^{n_i} \left( \frac{d^2}{64}\right)^{\alpha_{j,i}}\leq \left( \frac{d^2}{64}\right)^{\ld(\mu(D_i))}= \mu(D_i)^{2\ld(d)-6}\leq \mu(C)^ {2\ld(d)-6}.
$$

Furthermore, we have
$$
\prod_{j=1}^{n_i} \left(p_j^{\rho\ld(p_j)+\epsilon} \right)^{\alpha_{j,i}} \leq \left(\prod_{j=1}^{n_i}  
p_j^{\alpha_{j,i}} \right)^{\rho\ld(p_{n_i})+\epsilon}= \mu(D_i)^ {\rho\ld(p_{n_i})+\epsilon}\leq 
%\mu(C)^ {\rho\ld(p_{n_i})+\epsilon}\leq 
\mu(C)^{\rho \tau\ld(e)d+\epsilon},
$$
where the last inequality follows from $p_j <e^{\tau d}$ for all $j$ and equation \eqref{UpBoundMu}. 

Putting it all together, we get that
$$
\Hilb(F_k^i) \subset \left(\frac{d}{4} \cdot \mu(C)^{\gamma d+2\ld(d)+\kappa}\right)\Delta_C,
$$
where $\gamma= \rho\tau\ld(e) \approx 0.53$ and $\kappa= \epsilon -5 \approx 0.88$.
\end{proof}

Via simplification of the above notation we finally get

\begin{corollary}\label{FinResII}
Every simplicial $d$-cone
$C = \RR_+ v_1 + \cdots + \RR_+ v_d \subset \RR^d , \, d\geq 2$, which is not already unimodular (i.e., $\mu(C)>1$) has a unimodular triangulation $C = D_1 \cup \ldots \cup D_t$ such that for all $i$
$$
\Hilb(D_i) \subset \mu(C)^{\fFunc(d)}\Delta_C
$$
with $\fFunc(d) \in \bigO(d)$.
\end{corollary}

\section*{Acknowledgement} I thank the two anonymous referees for their careful reading of the paper. It led to improvements in the exposition, and helped me to correct an error in Lemma \ref{PrimeTransfer}.

\end{document}